\documentclass[12pt,amsart]{article}
\usepackage{amsmath, amssymb, amsfonts, amsthm,amscd}

\swapnumbers
\theoremstyle{plain}
\newtheorem{theorem}{Theorem}[section]
\newtheorem{lemma}[theorem]{Lemma}
\newtheorem{Titletheo}[theorem]{}

\theoremstyle{definition}
{\newtheorem{remark}[theorem]{Remark}

\numberwithin{equation}{theorem}
\renewcommand{\thetheorem}{\arabic{section}.\arabic{theorem}}

  \newcommand{\om}{\omega}    \newcommand{\Om}{\Omega}
  \newcommand{\sig}{\sigma}   
  \newcommand{\al}{\alpha}
  \newcommand{\del}{\delta}  
  \newcommand{\gam}{\gamma}   
  \newcommand{\lam}{\lambda}   
  \newcommand{\eps}{\varepsilon}
   \def\part{\partial}
  \newcommand{\wed}{\wedge}
  \newcommand{\bsym}{\boldsymbol}

  \renewcommand{\tilde}{\widetilde}
  \def\avi{\alpha^{\vee}_i}
  
  \def\uv{^{\vee}}

  \def\b1{\text{\bf\large 1}}  %big bold 1

  \def\simto{\overset{\sim}{\longrightarrow}}
  
  \def\ip<#1>{\langle#1\rangle}   %inner product--write `\ip<contents>'

  \newcommand{\Aff}{\operatorname{Aff}}

  \newcommand{\Imo}{\operatorname{Im}}
  
  \newcommand{\Ker}{\operatorname{Ker}}

\newcommand{\beqn}{\begin{equation}}
\newcommand{\eeqn}{\end{equation}}
\newcommand{\baln}{\begin{align}}
\newcommand{\ealn}{\end{align}}

\newcommand{\x}{\times}

 \newcommand{\fb}{\mathfrak{b}}
 \newcommand{\fg}{\mathfrak{g}}
 \newcommand{\fh}{\mathfrak{h}}
 \newcommand{\fn}{\mathfrak{n}}

 \newcommand{\fu}{\mathfrak{u}}

\newcommand{\bc}{\mathbb{C}}

\newcommand{\br}{\mathbb{R}}
\newcommand{\bz}{\mathbb{Z}}

 \newcommand{\cb}{\mathcal{B}}

 \newcommand{\cg}{\mathcal{G}}
 \newcommand{\ch}{\mathcal{H}}
 \newcommand{\ci}{\mathcal{I}}

 \def\cp{\wp}   %{\mathcal{P}}

 \newcommand{\ct}{\mathcal{T}}
 \newcommand{\cy}{\mathcal{Y}}
 \newcommand{\cz}{\mathcal{Z}}

\begin{document}

  \title{On  Cachazo-Douglas-Seiberg-Witten Conjecture for Simple Lie 
Algebras}
  \author{Shrawan Kumar\\Department of Mathematics\\
University of North Carolina\\
Chapel Hill, NC  27599--3250}
  \maketitle

\section{Introduction}

Let $\fg$ be a finite dimensional simple Lie algebra over the complex 
numbers $\bc$.  Consider the exterior algebra $R := \wed (\fg\oplus\fg )$ on 
two copies of $\fg$.  Then, the algebra $R$ is bigraded with the two copies 
of $\fg$ sitting in bidegrees (1,0) and (0,1) respectively.  To 
distinguish, we will denote the first copy of $\fg$ by $\fg_1$ and the 
second copy of $\fg$ by $\fg_2$.  Let $R^{p,q}$ denote the subspace of $R$ 
consisting of elements of bidegree $(p,q)$.  Thus,
  \[
R^{p,q} = \wed^p(\fg_1) \otimes \wed^q(\fg_2)\quad\text{ and }\quad
R = \oplus_{p,q\in\bz_+}\,R^{p,q}.
  \]

The bigrading, of course, gives rise to a $\bz_+$-grading by declaring any 
element of $R^{p,q}$ to have {\em total degree} $p+q$.  Set
  \[
R^n := \oplus_{p+q=n} \,R^{p,q};
  \]
thus $R^n$ consists of elements of total degree $n$.

The diagonal adjoint action of $\fg$ gives rise to a $\fg$-algebra 
structure on $R$ compatible with the bigrading.  We isolate three 
`standard' copies of the adjoint representation $\fg$ in $R^2$.  The 
$\fg$-module map 
  \[
  \part : \fg \to \wed^2(\fg ), \qquad x\mapsto \part x = \sum_i [x,e_i]\wed 
f_i, 
  \]
considered as a map to $\wed^2(\fg_1)$ will be denoted by $c_1$, and 
similarly, 
  \begin{align*}
c_2: &\fg \to \wed^2(\fg_2), \,\,\text{and}\\
c_3 : & \fg \to \fg_1\otimes\fg_2, \quad x\mapsto \sum_i\,
[x,e_i]\otimes f_i,
  \end{align*}
where $\{ e_i\}_{i\leq i\leq N}$ is any basis of $\fg$ and $\{ 
f_i\}_{1\leq i\leq N}$ is the dual basis of $\fg$ with respect to the 
normalized Killing form $\ip<\; ,\; >$ of $\fg$ (normalized as below).  It is easy to see that the three 
embeddings $c_i$ do not depend upon the choice of the basis $\{ e_i\}$.  
We fix a Cartan subalgebra $\fh$ of $\fg$ and a Borel subalgebra $\fb 
\supset \fh$, and normalize the Killing form by demanding that 
$\ip<\theta ,\theta > = 2$, for the highest root $\theta\in\fh^*$.
We denote by $C_i$ the image of $c_i$. 

Let $J$ be the (bigraded) ideal of $R$ generated by the three copies $C_1, 
C_2, C_3$ of $\fg$ (in $R^2$) and define the bigraded $\fg$-algebra
  \[
A := R/J.
  \]
The Killing form gives rise to a $\fg$-invariant $S\in A^{1,1}$ given by
  \[
S := \sum_i e_i\otimes f_i.
  \]

Motivated by supersymmetric gauge theory, 
Cachazo-Douglas-Seiberg-Witten [CDSW] made the following conjecture.  They 
proved the conjecture in  [CDSW], [W] for classical $\fg$.  Subsequently, 
Etingof-Kac proved the conjecture for $\fg$ of type $G_2$ by using the 
theory of abelian ideals in $\fb$.

Recall that the {\it dual Coxeter number} $h$ of $\fg$ is, by definition, 
$\langle\rho,\theta^\vee\rangle+1$, where $\rho$ is half the sum of all positive roots and $\theta^\vee$ is the  coroot corresponding to the highest root
$\theta$. The value of $h$  for any simple $\fg$ is given as follows (in the bracket): $A_{\ell} (\ell+1); B_{\ell} (2\ell-1); C_{\ell} (\ell+1);
D_{\ell} 
 (2\ell-2); E_6 (12); E_7 (18); E_8 (30); G_2 (4); $ and   
 $F_4 (9).$

  \begin{Titletheo}{\rm\bf Conjecture [CDSW]} (i) The subalgebra $A^{\fg}$ of
$\fg$-invariants in $A$ is generated, as an algebra, by the element $S$.

(ii) $S^h =0$.

(iii) $S^{h-1}\neq 0$. 

Thus, as an algebra,
  \[
A^{\fg} \simeq \bc [S]/\ip<S^h>,
  \]
where $\ip<S^h>$ denotes the ideal of the polynomial ring $ \bc [S]$ generated
by $S^h$. 
\end{Titletheo}

The aim of this paper is to give a uniform proof of the above conjecture part (i). In addition, we give a conjecture (cf. Section 3), the validity of which would imply  
part (ii) of the above conjecture. 

The proof of part (i) follows  from our Theorem 2.2, which is the main theorem of this paper. This theorem identifies the graded algebra $B^\fg$ with the singular cohomology of a certain (finite dimensional) projective 
subvariety  $\cy_2$ of the infinite Grassmannian   $\cy$ associated to $\fg$, where 
$B := R/\ip<C_1\oplus C_2>$. The definition of the subvariety  $\cy_2$
is motivated from the theory of abelian ideals in the Borel subalgebra $\fb$
of $\fg$. 
This theorem is proved by using Garland's result 
on the Lie algebra cohomology of  
$\hat{\fu} := \fg\otimes t\bc [t]$; Kostant's result on the `diagonal' 
cohomolgy of $\hat{\fu}$ and its connection with abelian ideals in $\fb$;
and a certain deformation of the singular cohomology of  $\cy$ introduced by
Belkale-Kumar [BK]. 

To obtain part (i) of the above conjecture, observe that the singular cohomology $H^*(\cy)$ surjects onto $H^*(\cy_2)$ under the restriction map. Moreover, as is well known, $H^*(\cy)$ (being isomorphic with the cohomology of the based 
loop space $\Omega_e(K)$ of a maximal compact subgroup $K$ of $G$) is  isomorphic with the polynomial ring $\bc[x_1,\ldots, x_\ell]$ in $\ell$ variables, where $G$ is the connected, simply-connected complex algebraic group with 
Lie algebra $\fg$ and $\ell$ is its rank. Thus, by virtue of our Theorem 2.2, we get a surjective algebra 
homomorphism from  $\bc[x_1,\ldots, x_\ell]$ onto $H^*(\cy_2)\simeq B^\fg$. Finally, we show that each $x_i$ (for $i> 1$) goes to zero under the canonical quotient map $B^\fg\to A^\fg$. This is proved by using the explicit expression of the suspension maps 
associated to the universal principal $K$-bundle $K\to E(K) \to B(K)$ and
also the fibration $\Om_e(K) \to P_e(K) \to K$, where $P_e(K)$ denotes the
space of continuous maps $\beta : I\to K$ from the closed unit interval
$I$ starting at $e$. Moreover, it is easy to see that  $x_1$ goes to the element $S$ (up to a nonzero scalar multiple). This proves  part (i) of the above conjecture.
\vskip2ex

It is my pleasure to thank P. Etingof for bringing to my attention the above conjecture and explaining to me his work with Kac. I also thank B. Kostant and R. Suter for some helpful correspondences. This work was partially supported by the NSF grant DMS 0401084. 

  \section{A topological model for $\Bigl( \frac{R}{\ip<C_1\oplus 
C_2>}\Bigr)^{\fg}$}

\subsection{Notation and Preliminaries}

Let $G$ be the connected, simply-connected complex algebraic group with 
Lie algebra $\fg$ and let $\cg$ by the corresponding affine Kac-Moody 
group which is, by definition, a $\bc^*$-central extension of the loop 
group $L(G) := G\Bigl( \bc [[t]][t^{-1}]\Bigr)$, where  $G\Bigl( \bc [[t]][t^{-1}]\Bigr)$ denotes the $ \bc [[t]][t^{-1}]$-rational points of the algebraic group $G$ (cf. [K2, Section 13.2]).  
Let $\cp$ be the standard maximal parabolic subgroup, which is by 
definition the $\bc^*$-extension of $G\bigl( \bc [[t]]\bigr)$.  Then, the 
{\em infinite Grassmannian} $\cy$ is, by definition,
  \[
\cy := \cg /\cp ;
  \]
this is equipped with a projective ind-variety. The ring homomorphism
$e:\bc[[t]] \to \bc, Q \mapsto Q(0),$ gives rise to a group homomorphism
$\hat{e}:\cp \to G$.  Let $\cb$ be the standard 
Borel subgroup of $\cg$ defined as $\hat{e}^{-1}(B)$, $B$ being the Borel subgroup of $G$ with Lie algebra $\fb$.  By the Bruhat decomposition 
(cf. [K2, $\S$13.2.12]), $\cy$ is the 
disjoint union of the Bruhat cells
  \[
\cy = \bigsqcup_{w\in\Aff W/W} \cb w \cp /\cp ,
  \]
where $W$ is the Weyl group of $G$ and $\Aff W$ is the corresponding affine Weyl group, 
which is by definition the semidirect product of $W$ with the coroot 
lattice $Q\uv\subset\fh$ acted on via the natural action of $W$ on $\fh$.  
$\Aff W$ acts on $\fh$ via the standard action of $W$ on $\fh$ and $Q\uv$ 
acts on $\fh$ via translation.

Let $\Aff' (W)\subset \Aff (W)$ be the set of minimal coset 
representatives of the $W$-cosets  $\Aff W/W$.  Define the Weyl 
chamber
  \[
\fh_+ := \{ h\in\fh : \al_i(h)\in\br_+ \text{ for all simple roots 
}\al_i\} .
  \]
Let $C := \{ h\in\fh_+ :  
\theta (h)\leq 1 \}$ be the fundamental alcove.  Then, $C$ is a 
fundamental domain for the action of $\Aff (W)$ on $\fh$.  

Moreover, as is well known,
  \[
\Aff' (W) = \{ w\in\Aff (W): w^{-1}C \subset \fh_+\} .
  \]
We define a closed (finite dimensional) subvariety $\cy_2$ of $\cy$ by
  \[
\cy_2 := \bigcup_{w\in\Aff_2'(W)}\,  \cb w \cp /\cp ,
  \]
where
  \[
\Aff'_2(W) := \{ w\in\Aff (W) : w^{-1}C \subset 2C\} .
  \]
Let $\ell$ be the rank of $\fg$, i.e., the dimension of the Cartan subalgebra
$\fh$. Since volume $(2C) = 2^{\ell}\; \text{volume}\, (C)$,
  \beqn
\# \Aff'_2(W) = 2^{\ell} .
  \eeqn
In general, $\cy_2$ is not an irreducible variety.

For any topological space $X$, we denote by $H^*(X)$ the singular cohomology 
of $X$ with complex coefficients. 

Following is the main result of this paper.
 \stepcounter{theorem}

  \begin{theorem} The singular cohomology $H^*(\cy_2)$ of $\cy_2$
 is isomorphic as a $\bz_+$-graded algebra with
the graded algebra $B^\fg$, where $B := R/\ip<C_1\oplus C_2>$ and
 $\ip<C_1\oplus C_2>$ denotes the ideal of $R$
generated by $C_1\oplus C_2$.
  \end{theorem}

Before we come to the proof of the above theorem, we need to define a new product $\odot_0$ in the singular cohomology $H^*(\cy)$
of $\cy$ as follows.  (This was introduced in [BK] for  any flag variety of any
semisimple group $\ch$ 
in an attempt to solve the `Hermitian eigenvalue problem' for $\ch$.)

Let $\{\eps^u\}_{u\in\Aff'(W)}$ be the Schubert basis of $H^*(\cy)$
defined by $\eps^u(y_v)=\delta_{u,v},$ for any ${v\in\Aff'(W)}$, where 
$y_v$ is the image in the homology $H_*(\cy)$ of the fundamental class of the closure $\overline{ \cb v \cp /\cp}$.  
Express  the
usual cup product
  \beqn 
\eps^u\cdot \eps^v = \sum_w c^w_{u,v} \eps^w \,.
  \eeqn

Let $\hat{\fg}$ be the affine Kac--Moody Lie algebra associated to $\fg$. 
Recall that 
$$\hat{\fg}:= \fg\otimes \bc[t,t^{-1}] \oplus \bc c \oplus \bc d$$
with the bracket
defined by
\begin{align*}  
\bigl[ t^m\otimes x+\lam c+\mu d, t^{m'}\otimes x' +\lam'c
+\mu'd\bigr] &= \\
\bigl( t^{m+m'}\otimes [x,x'] + \mu m' t^{m'}\otimes x' - &\mu'm
t^m\otimes x\bigr) \oplus m\del_{m,-m'} \langle x,x'\rangle c, 
\end{align*}
for $\lam ,\mu ,\lam', \mu'\in\bc\,, \,
m,m'\in\bz$ and $x,x'\in\fg$.  
  
Let $\ct\subset\cg$ be the standard maximal torus of $\cg$ and let 
$\hat{\fh}$ be its Lie algebra.  Then, $\hat{\fh}= \fh\otimes 1\oplus 
\bc c \oplus \bc d$. Let $\{\al_0,\al_1, \ldots, \al_{\ell} 
\}$ be the simple roots of $\hat{\fg}$ such that  $\{\al_1, \ldots, \al_{\ell} 
\}$ are the simple roots of $\fg$.  Thus, $\al_i\in (\hat{\fh})^*$.  Choose 
an element $x_j\in\hat{\fh}$, $0\leq j\leq\ell$, satisfying $\al_i(x_j) = 
\del_{i,j}$, for all $0\leq i\leq\ell$.  Also, let $\hat{\rho}\in 
(\hat{\fh})^*$ be an element satisfying $\hat{\rho}(\avi)=1$ for all 
simple coroots $\avi$, $0\leq i\leq\ell$.  Now, the new product $\odot_0$ 
is defined by
  \[
\eps^u\odot_0\eps^v = \sum_{w\in\Aff' (W)} c^w_{u,v} \del_{d^w_{u,v},0} 
\,\eps^w, 
  \]
where $\del$ is the Kronecker delta and 
  \[
d^w_{u,v} := \bigl( u^{-1}\hat{\rho}+v^{-1}\hat{\rho} -w^{-1}\hat{\rho} - 
\hat{\rho}\bigr) (x_0).
  \]
As proved in [BK], this product is associative (and of course commutative).

In fact, one can define this new product $\odot_0$ in the cohomology of 
$\ch /\ci$ for any Kac-Moody group $\ch$ and any standard parabolic 
subgroup $\ci$ and for the cohomology with integral coefficients $\bz$. However, we need this 
product only for $\cy$. Moreover, there is a $q$-deformation of the cup product in $H^*(\ch /\ci, \bz)$ such that the usual cup product corresponds to the value $q=1$, whereas the new product  $\odot_0$ corresponds to the value $q=0$
(cf. [BK, $\S$6]).  

For any left $\cb$-stable closed subset $\cz\subset\cy$, the kernel of the
restriction map $\gam : H^*(\cy ) \to H^*(\cz )$ is an ideal under the
product $\odot_0$ in $H^*(\cy )$.  To see this, let $\cz =\bigcup \cb
w\cp$, where $w$ runs over a certain subset $S_{\cz}$ of $\Aff' (W)$.  
Then,
  \beqn
\Ker \gam = \bigoplus_{w\in\Aff' (W)\backslash S_{\cz}} \bc\,\eps^w . 
\tag{$*$}
  \eeqn

Of course, $\Ker \gam$ is an ideal of $H^*(\cy )$ under the cup product.  
From this and $(*)$, it follows easily that $\Ker \gam$ is an ideal of 
$H^*(\cy )$ under $\odot_0$.

This allows us to define the product $\odot_0$ in $H^*(\cz )$ making the 
restriction map $\gam : H^*(\cy )\to H^*(\cz )$ a ring homomorphism under 
$\odot_0$.  In particular, we have the product $\odot_0$ in $H^*(\cy_2)$.

  \begin{lemma}
The product $\odot_0$ in $H^*(\cy_2)$ coincides with the cup product.
  \end{lemma}

  \begin{proof}  From the definition of $\odot_0$, to prove the lemma, it 
suffices to show that for any $u,v,w\in\Aff'_2(W)$ such that 
$c^w_{u,v}\neq 0$ in the decomposition (1) of (2.2), we have
  \beqn
\bigl(u^{-1}\hat{\rho} + v^{-1}\hat{\rho} - w^{-1}\hat{\rho} -\hat{\rho}\bigr)(x_0) = 
0.
  \eeqn
Since $c^w_{u,v}\neq 0$, we, in particular, have
  \beqn
\ell (w) = \ell (u) +\ell (v).
  \eeqn
By [Ko2, Theorems 3.13 and 4.5 and Identity 3.18], for any 
$u\in\Aff'_2(W)$, 
  \beqn
\rho - u^{-1}\rho = \ell (u)\del + \lam^u ,
  \eeqn
for some $\lam^u$ in the root lattice of $\fg$, where $\del$ is the basic 
imaginary root of the affine Lie algebra $\hat{\fg}$.

Thus, for any such $u$, 
  \beqn
(\rho - u^{-1}\rho )(x_0) = \ell (u).
  \eeqn

Combining (2) and (4), we get (1).  This proves the lemma.
  \end{proof}

We now come to the proof of Theorem 2.2.

\stepcounter{theorem}
  \begin{proof}[{\rm\bf 2.4}\quad Proof of Theorem 2.2] Define the following
Lie 
subalgebras of $\hat{\fg}$:
  \[
\hat{\fu} := \fg\otimes t\bc [t], \quad \hat{\fu}^- = \fg \otimes t^{-1} 
\bc [t^{-1}], \quad \fg_0:=\fg\otimes 1\oplus \bc c \oplus \bc d .
  \]

By a result of Garland (cf. [K2, Theorem 3.2.7]), as 
$\fg_0$-module, the Lie algebra cohomology
  \beqn
H^p(\hat{\fu}^-) \simeq \bigoplus_{w\in\Aff'(W), \ell (w)=p} 
L_0(w^{-1}\hat{\rho}-\hat{\rho})^*, 
  \eeqn
where $L_0(\lam )$ denotes the irreducible $\fg_0$-module with highest 
weight $\lam$.  Similarly, by [K2,(3.2.11.3)], 
  \beqn
H^p(\hat{\fu}) \simeq \bigoplus_{w\in\Aff' (W), \ell (w)=p} 
L_0(w^{-1}\hat{\rho}-\hat{\rho}).
  \eeqn

By [BK, Theorem 42], there is a  graded algebra isomorphism
  \[
\phi : H^*(\cy, \odot_0) \simeq \bigl[ H^*(\hat{\fu})\otimes 
H^*(\hat{\fu}^-)\bigr]^{\fg_0}, 
  \]
where  $H^*(\cy, \odot_0)$ denotes   $H^*(\cy)$ equipped with  the product 
$\odot_0$.  (Even though in [BK] the 
theorem is proved for semisimple Lie algebras, the same proof carries 
over and gives the result for any symmetrizable Kac-Moody case.)  Moreover, 
under the above identifications (1) and (2), for any $w\in\Aff' (W)$,
  \[
\phi (\eps^w) \in \bigl[ L_0(w^{-1}\hat{\rho}-\hat{\rho}) \otimes 
L_0(w^{-1}\hat{\rho}-\hat{\rho})^*\bigr]^{\fg_0} .
  \]
The cohomology modules $H^p(\hat{\fu})$ and $H^p(\hat{\fu}^-)$ acquire a 
grading coming from the total degree of $t$ in $\wed^p(\hat{\fu})$ and 
$\wed^p(\hat{\fu}^-)$ respectively.  Thus, 
  \[
H^p(\hat{\fu}) = \bigoplus_{m\in\bz_+} H^p_{(-m)} (\hat{\fu}) ,
  \]
where $H^p_{(-m)}(\hat{\fu})$ denotes the space of elements of 
$H^p(\hat{\fu})$ of total $t$-degree $-m$.  Similarly,
  \[
H^p(\hat{\fu}^-) = \bigoplus_{m\in\bz_+} H^p_{(m)} (\hat{\fu}^-).
  \]
Clearly, 
  \[
H^p_{(-m)}(\hat{\fu}) =  H^p_{(m)}(\hat{\fu}^-) =0, \quad\text{ if } m<p.
  \]
Let $H^*_D(\hat{\fu})$ denote the `diagonal' cohomology
  \[
H^*_D(\hat{\fu}) := \bigoplus_{p\in\bz_+} H^p_{(-p)}(\hat{\fu}), 
  \]
and similarly
  \[
H^*_D(\hat{\fu}^-) := \bigoplus_{p\in\bz_+} H^p_{(p)}(\hat{\fu}^-).
  \]
Then, clearly, $H^*_D(\hat{\fu})$ and $H^*_D(\hat{\fu}^-)$ are graded $\fg$-subalgebras of  $H^*(\hat{\fu})$ and $H^*(\hat{\fu}^-)$ respectively. By [Ko2, Theorem 4.16], as graded $\fg$-algebras, 
   \[
H^*_D(\hat{\fu}) \simeq \wed (\fg ) \big/ \ip<\part (\fg )> ,
  \]
where $\part : \fg\to\wed^2(\fg )$ is the map defined in Section (1) and 
$\ip<\part (\fg )>$ denotes the ideal of the exterior algebra $\wed (\fg 
)$ generated by $\part (\fg )$.  Similarly, 
  \[
H^*_D (\hat{\fu}^-) \simeq \frac{\wed (\fg )}{\ip<\part (\fg )>} .
  \]
Moreover, by [Ko1], the $\fg$-module $\frac{\wed (\fg )}{\ip<\part (\fg 
)>}$
is a multiplicity free $\fg$-module.  In fact, for any abelian ideal 
$I\subset \fb$ (which is automatically contained in $\fn := [\fb ,\fb ]$), consider the line  $\wed^{\dim I}(I) 
\subset \wed^{\dim I}(\fg )$.  Then, $I$ being an ideal in 
$\fb$, $\wed^{\dim I}(I)$ is stable under $\fb$.  Let $V_I\subset 
\wed^{\dim I}(\fg )$ be the (irreducible) $\fg$-module generated by 
$\wed^{\dim I}(I)$.  Then, as $\fg$-modules, the canonical map (induced by the
inclusion)
  \[
\bigoplus _{I\in\Xi} V_I \to \frac{\wed (\fg )}{\ip<\part (\fg )>} 
  \]
is an isomorphism, where $\Xi$ is the set of all abelian ideals of $\fb$.  By a result of D. Peterson, $\Xi$ has cardinality 
$2^{\ell}$, $\ell$ being the rank of $\fg$.  Further,  for any 
$p\geq 0$, 
  \[
\frac{\wed^p (\fg )}{\ip<\part (\fg )>\cap \wed^p(\fg )} \;\text{ is a 
self-dual $\fg$-module}.
  \]
Thus, as graded algebras,
  \beqn
\bigl[ H^*_D (\hat{\fu}) \otimes H^*_D(\hat{\fu}^-)\bigr]^{\fg_0} \simeq 
\Biggl[ \frac{\wed (\fg )}{\ip<\part \fg >} \otimes \frac{\wed (\fg 
)}{\ip<\part \fg >}\Biggr]^{\fg}=:B^\fg.
  \eeqn

We next show that the composite $\xi := \gam\circ\phi^{-1}\circ i$ 
  \begin{align*}
\bigl[ H^*_D(\hat{\fu})\otimes H^*_D(\hat{\fu}^-)\bigr]^{\fg_0} 
&\overset{i}{\hookrightarrow} \bigl[ H^*(\hat{\fu}) \otimes 
H^*(\hat{\fu}^-)\bigr]^{\fg_0} \\
&\overset{\phi^{-1}}{\simto} H^*(\cy ,\odot_0) 
\overset{\gam}{\longrightarrow} H^*(\cy_2,\odot_0),
  \end{align*}
is an algebra isomorphism.
   
By (1) of (2.1), the total dimension 
  \[
\dim H^*(\cy_2) = 2^{\ell}.
  \]
Also, by the above discussion, 
  \[
\dim \bigl[ H^*_D(\hat{\fu}) \otimes H^*_D(\hat{\fu}^-)\bigr]^{\fg_0} = 
2^{\ell}.
  \]
Observe that all the maps $i$, $\phi^{-1}$ and $\gam$ are algebra 
homomorphisms.  Thus, it suffices to show that $\xi$ is surjective.  Take 
$w\in\Aff'_2(W)$ and consider the one dimensional subspace
  \[
V(w) := \bigl[ L_0(w^{-1}\hat{\rho}-\hat{\rho}) \otimes 
L_0(w^{-1}\hat{\rho} - \hat{\rho})^*\bigr]^{\fg_0} \text{  of  } \bigl[ 
H^*(\hat{\fu}) \otimes H^*(\hat{\fu}^-)\bigr]^{\fg_0}
  \]
under the identifications (1) and (2) of (2.4).  Under the isomorphism 
$\phi^{-1}$, 
  \[
\phi^{-1}(V(w)) = \bc \,\eps^w .
  \]
On the other hand, by [Ko2, Theorems 4.5 and 4.8], for $w\in\Aff'_2(W)$, 
$V(w)\subset \Imo(i).$  This proves the surjectivity of $\xi$ and hence concludes the proof that $\xi$ is a graded algebra isomorphism. Combining this with the identification (3) and Lemma 2.3, we get Theorem 2.2. 
\end{proof}

 By a result of Garland-Raghunathan [GR] and Quillen (unpublished), the space $\cy$ is homotopically equivalent to the continuous based loop group $\Om_e(K)$ of a maximal compact subgroup $K$ of $G$.  In particular, 
  \[
H^*(\cy) \simeq H^*(\Om_e(K)).
  \]
 \stepcounter{theorem}
 \renewcommand{\thesubsection}{\thetheorem}

  \subsection{The generators of {\bf $H^*(\bsym{\Om}_e(K))$}} Let $F$ be
the ring $\bc [t,t^{-1}]$ of Laurent polynomials.  For any invariant
homogeneous polynomial $P\in S(\fg^*)^\fg$ of degree $d+1$ ($d\geq 1$), define the map
  \[
\tilde{\phi}_P : \wedge^{2d}_\bc (\fg\otimes F) \to \Om^1_F
  \]
by
  \begin{multline*}
\tilde{\phi}_P \Bigl( v_0\wed v_1\wed\cdots\wed v_{2d-1}\Bigr) =\\
 \sum_{\sig\in S_{2d}} \eps (\sig )\, P\Bigl( v_{\sig (0)}, \bigl[ v_{\sig
(1)}, v_{\sig (2)}\bigr] , \ldots, \bigl[ v_{\sig(2d-3)}, v_{\sig (2d-2)}
\bigr] , dv_{\sig (2d-1)}\Bigr) ,
  \end{multline*} for $v_i\in\fg\otimes F$, where the loop algebra 
$\fg\otimes F$ has the standard bracket,  $\Om^1_F$ is the space of
algebraic 1-forms on the affine variety $\bc^*$, $d(x\otimes a(t))$ $:=
x\otimes \frac{da(t)}{dt}\, dt\in\fg\otimes\Om^1_F$ for $x\in\fg$ and
$a(t)\in F$, and $P$ is extended $F$-linearly.  Define the residue map
$\al : \Om^1_F \to\bc$ by
  \[ \al \bigl( Q(t)dt\bigr) = \text{the coefficient of $\frac{1}{t}$ in
$Q(t)$}.  \] Of course, $\al (Q(t)dt) = \frac{1}{2\pi i} \int_{S^1}
Q(t)dt$.  Composing $\tilde{\phi}_P$ with $\al$, we get the map
  \[
\phi_P : {\wedge}_{\bc}^{2d} (\fg\otimes F) \to \bc .
  \]

It is easy to see that $\phi_P$ represents a relative cocycle for the Lie
algebra pair $(\fg\otimes F, \fg )$ (cf. [K2, \S3.1.3]), thus giving rise
to a Lie algebra cohomology class $[\phi_P ]\in H^{2d}(\fg\otimes F, \fg
)$.

We recall the following result from [K1, Theorem 1.6 and Lemma 1.8].

  \begin{theorem} The integration map (appropriately defined) defines an
algebra isomorphism in cohomology
  \[
H^*(\fg\otimes F, \fg ) \simeq H^*(\cy ) \simeq H^*
\bigl(\Om_e(K)\bigr). 
  \]
 \end{theorem}

From now on, we will identify $H^*(\fg\otimes F,\fg )$ with $H^*(\cy )$.

As is well known, $S(\fg^*)^\fg$ is freely generated (as a commutative
algebra) by certain homogeneous polynomials $P_1, \ldots, P_{\ell}$ of
degrees $m_1+1, m_2+1,\ldots, m_{\ell}+1$ respectively, where $m_1=1 <
m_2 \leq \cdots \leq m_{\ell}$ are the exponents of $\fg$.  Moreover,
$H^*(\Om_e(K))$ is freely generated (again as a commutative algebra) by
taking double suspension $\hat{P}_1,\ldots, \hat{P}_{\ell}$
of the elements $P_1,\ldots, P_{\ell}$ respectively; first
associated to the universal principal $K$-bundle $K\to E(K) \to B(K)$ and
then to the fibration $\Om_e(K) \to P_e(K) \to K$ where $P_e(K)$ denotes the
space of continuous maps $\beta : I\to K$ from the closed unit interval
$I$ starting at $e$ (i.e., $\beta (0)=e$) (cf. [B, Proposition 7.1]).
(Here we have identified $H^*(B(K))\simeq S(\fg^*)^\fg$.) Observe that
$\hat{P}_i\in H^{2m_i}(\Omega_e(K))$. 

Then, the class $\hat{P}_i$ corresponds to the class $[\phi_{P_i}]
\in  H^{2m_i}(\fg\otimes F, \fg )$ under the identification of Theorem 2.6 
 (cf. [CS, Proposition 3.2] where an explicit description of the first suspension is given, combined with [B, $\S$7] where the second suspension is 
described explicitly; alternatively see [T, $\S$3]). 

Combining this with Theorem 2.6, we get the following.

  \begin{theorem} The cohomology classes $[\phi_{P_1}], \ldots, [\phi_{P_{\ell}}] \in H^*(\fg\otimes F, \fg )$ (with degrees 
 $2m_1,\ldots, 2m_{\ell}$ respectively)  generate the cohomology $H^*(\fg\otimes F,
\fg )$ freely (as a commutative algebra).
  \end{theorem}

Theorems 2.2 and 2.7 readily prove the first part of Conjecture 1.1.  
Specifically, we have the following.

  \begin{theorem} With the notation as in Section (1), the algebra $A^{\fg}$
is generated, as an algebra, by the element $S$.
    \end{theorem}

  \begin{proof}  By the above Theorem 2.7, the polynomial ring 
  \[
\bc [x_1,\ldots, x_{\ell}] \simto H^*(\cy ), \; x_i \mapsto [\phi_{P_i}].
  \] Observe next that the restriction map $\gam : H^*(\cy )\to
H^*(\cy_2)$ is surjective.  This follows since, by the Bruhat
decomposition of $\cy$ as in (2.1), $\cy$ is a CW-complex obtained from $\cy_2$ by
attaching (real) even dimensional cells.  In particular, composing the above isomorphism with $\gamma$, we get a
surjection
  \[
\eta : \bc [x_1,\ldots, x_{\ell}] \to H^*(\cy_2).
  \]

We next show that, under the identification of $H^*(\cy_2)$ with the
algebra $B^\fg$ (guaranteed by Theorem 2.2), the element $x_1 \mapsto zS$ (for
a nonzero $z\in\bc$) and each $\eta (x_i)$, $i>1$, lies in the kernel of the
standard quotient map $B^{\fg} \to A^\fg$.

We first prove the assertion that $\eta (x_1)=zS$.  Since $H^2(\cy_2)$ is
one dimensional and $\eta (x_i)$ is a cohomology class of degree $2m_i >2$
for any $i>1$, $\eta(x_1)$ is the unique nonzero element of 
$H^2(\cy_2)$ (up to a nonzero scalar multiple). Further, by Theorem 2.2,
  \[
H^2(\cy_2) \simeq (B^2)^{\fg},
  \] 
where $B^2$ is the total degree 2 component of $B$.  But $S$ is the
unique element of $(B^2)^\fg$ (up to scalar multiples).  This proves the
assertion that $\eta (x_1)=zS$.

We next prove that $\eta (x_i)$, for $i>1$, lies in the kernel of $B^\fg
\to A^\fg$.  
  \bigskip

Take dual bases $\{ e_i\}$ and $\{ f_i\}$ of $\fg$ as in Section (1).  Then,
$\{ e_i(n) := e_i\otimes t^n\}_{n\in\bz}$ forms a basis of $\fg\otimes F$.  
Identify $(\fg\otimes t^n)^*$ with $\fg\otimes t^n$ as a $\fg$-module via
the Killing form on the $\fg$-factor.  Then, $\{ f_i(n) := f_i\otimes
t^n\}_{n\in\bz}$ is the dual basis of
  \beqn (\fg\otimes F)\uv := \bigoplus_{n\in\bz} (\fg\otimes t^n)^* \simeq
\fg\otimes F .
  \eeqn

For any $\sig\in S_{2d}$ and $P\in S^{d+1}(\fg^*)^\fg$ ($d\geq 2$), consider the 
linear form
  \[
\phi_{P,\sig}: \otimes^{2d}_\bc (\fg\otimes F) \to \bc ,
  \]
defined by
 \begin{multline*} {\phi}_{P,\sig} \Bigl( v_0\otimes v_1\otimes
\cdots\otimes v_{2d-1}\Bigr)\\
  = \alpha\Bigl(P\Bigl( v_{\sig (0)}, \bigl[ v_{\sig 
(1)}, v_{\sig (2)}\bigr] , \ldots, \bigl[ v_{\sig(2d-3)}, v_{\sig (2d-2)}
\bigr] , dv_{\sig (2d-1)}\Bigr)\Bigr).
  \end{multline*}
For notational convenience, assume $\sig (1) < \sig (2)$.  For any fixed 
$$v_0,v_1, \ldots, \hat{v}_{\sig (1)}, \ldots, \hat{v}_{\sig (2)}, 
\ldots, v_{2d-1} \in \fg\otimes F,$$
 consider the restriction 
$\bar{\phi}_{P,\sig}$ of the function $\phi_{P,\sig}$ to 
$$v_0\x v_1\x 
\cdots\x \fg\otimes F\x\cdots\x \fg\otimes F\x\cdots\x v_{2d-1},$$
 where 
the two copies of $\fg\otimes F$ are placed in the $\sig (1)$ and $\sig 
(2)$-th slots.  Then, under the above identification (1), 
  \begin{align*}
\bar{\phi}_{P,\sig} &= \sum_{i,j,m,n} f_i(n) \otimes f_j(m)\, \alpha\Bigl(P\Bigl( 
v_{\sig (0)}, \bigl[ e_i(n),e_j(m)\bigr] ,\bigl[ v_{\sig 
(3)}, v_{\sig (4)}\bigr] ,\ldots,\\
&\hspace{1.8in}\bigl[ v_{\sig (2d-3)}, v_{\sig (2d-2)} \bigr] , dv_{\sig 
(2d-1)}\Bigr)\Bigr) \\
  &= \sum_{i,j,m,n,k} f_i(n)\otimes f_j(m)\, \alpha\Bigl(P\bigl( -, \ip<[e_i,e_j], 
e_k> f_k{(n+m)}, -\bigr)\Bigr) \\
 &= \sum_{i,j,m,n,k} \ip<e_i, [e_j,e_k]>\, f_i(n)\otimes f_j(m)\, \alpha\Bigl(P\bigl( 
-, f_k{(n+m)}, -\bigr)\Bigr) \\
 &= \sum_{j,k,m,n} [e_j,e_k](n) \otimes f_j(m)\, \alpha\Bigl(P\bigl( -, f_k{(n+m)}, 
-\bigr)\Bigr)\\
 &= - \sum_{j,k,m,n} [e_k,e_j](n) \otimes f_j(m)\, \alpha\Bigl(P\bigl( -, f_k{(n+m)}, 
-\bigr)\Bigr) .  \tag{2}
  \end{align*}
 \stepcounter{equation}

Recall the definition of the isomorphism
  \[
\xi : \bigl[ H^*_D(\hat{\fu}) \otimes H^*_D(\hat{\fu}^-)\bigr]^{\fg_0} \to 
H^*(\cy_2)
  \]
from (2.4).  In particular, any class $c\in H^{2p}(\cy_2)$ is the 
restriction of a class in $H^{2p}(\cy )$ represented by a relative cocycle 
$\om$ of the Lie algebra pair $(\fg\otimes F,\fg )$ lying in the linear 
span of 
  \[
\Bigl\{ f_{i_1}(1)\wed \cdots \wed f_{i_p}(1)\wed f_{j_1}(-1) 
\wed\cdots\wed f_{j_p}(-1)\Bigr\}_{\begin{matrix} i_1<\cdots 
<i_p\\[-5pt] j_1 < \cdots <j_p \end{matrix}} .
  \]

Also, recall the isomorphism
  \[
\bigl[ H^*_D(\hat{\fu})\otimes H^*_D(\hat{\fu}^-)\bigr]^{\fg_0} \simeq 
B^\fg
  \]
from (3) of (2.4).  These two isomorphisms put together give rise to the 
identification of Theorem 2.2:
  \[
H^*(\cy_2) \simto B^\fg .
  \]

From the definition of the map $\eta : \bc [x_1,\ldots, x_{\ell}]$ $\to 
H^*(\cy_2)$ and (2), we see that, for all $i>1$, $\eta (x_i)$ belongs to 
the ideal of $B$ generated by the copy $C_3$ of the adjoint 
representation.  To prove this observe that, from the above discussion, in 
the decomposition (2) of $\bar{\phi}_{P,\sig}$ we can take $(n,m)$ to be 
one of (1,1), (-1,-1), (1,-1) or (-1,1).  The terms in the right side of 
(2) of the form (1,1) (resp. (-1,-1)) lie in the ideal generated by $C_1$ 
(resp. $C_2$), whereas the terms of the form (1,-1) and (-1,1) lie in the 
ideal generated by $C_3$.  Thus, $\eta (x_i)$ goes to 0 under the quotient 
map $B^\fg \to A^\fg$.  This proves the theorem.
  \end{proof}

\begin{remark} It is likely that 
 Etingof's map given in his conjecture [E, Conjecture 2.3] coincides with our map $\eta$ (defined in the proof of Theorem 2.8) under the identification of Theorem 2.2. If so, the surjectivity of the map $\eta$ as in the proof of  Theorem 2.8
together with  Theorem 2.2 will prove his conjecture.
\end{remark}

 \section{A Conjecture}

Let $\Xi$ be the set of all the abelian ideals in the Borel subalgebra $\fb$ and let $\Xi^o$ be the subset consisting of those abelian ideals
$I$ such that for any root space $\fg_\alpha$ corresponding to the root $\alpha$, if $\fg_\alpha \subset I$, then $\langle\alpha,\theta\rangle\neq 0$, where
$\theta$ is the highest root. 

Recall that there is a bijection  $ \zeta: \Xi \to \Aff'_2(W) 
 $ such that for any $I \in \Xi, \dim I =   \ell(\zeta (I))$
(see, e.g., [Ko2, Theorem 4.4]). Actually, we take  $\zeta (I)= (\sigma_I)^{-1}$, 
where  $\sigma$ is the map as in loc cit. Denote the image $\zeta (\Xi^o)$ by $ \Aff'_2(W)^o$ and consider the subvariety of $\cy_2$:
$$\cz:= \bigcup_{w\in\Aff_2'(W)^o}\,  \cb w \cp /\cp .$$

 \begin{Titletheo}{\rm\bf Conjecture.} Consider the map $\Phi: H^*(\cy_2) \to A^\fg$ obtained by the isomorphism $H^*(\cy_2) \simto B^\fg$ of Theorem 2.2 followed by the standard projection $B^\fg \to  A^\fg$. Then, $\Phi$ factors through the cohomology
$H^*(\cz)$ under the restriction map $H^*(\cy_2) \to H^*(\cz)$. 
\end{Titletheo} 

By Suter [S, Theorem 11 and Proposition 4], for any $I \in \Xi^o, \dim I \leq h-1$. In particular, this gives
$$\dim_\bc \cz\leq h-1.$$

So, the validity of the above conjecture will readily imply the validity of 
part (ii) of Conjecture 1.1.

  \end{document}